\documentclass[11pt]{article}
\usepackage{graphicx}
\usepackage{geometry}
\usepackage{tikz}
\usepackage{amsmath}

\newtheorem{theorem}{Theorem}
\newtheorem{observation}[theorem]{Observation}
\newtheorem{cor}[theorem]{Corollary}
\newtheorem{proposition}[theorem]{Proposition}
\newtheorem{conjecture}[theorem]{Conjecture}
\newtheorem{lemma}[theorem]{Lemma}
\newtheorem{question}[theorem]{Question}

\newenvironment{Proof}{\proofing}{\QED}
\newcommand{\QED}{\hspace{8mm}\mbox{\textsc{qed}}\medskip}
\newcommand{\proofing}{\textsc{Proof.}}

\newcommand{\FD}{\mathop{\mathit{FD}}}
\newcommand{\FTD}{\mathop{\mathit{FTD}}}
\newcommand{\FdT}{\mathop{\mathit{FDT}}}
\newcommand{\td}{\mathop{\mathit{td}}}
\newcommand{\dom}{\mathop{\mathit{dom}}}
\newcommand{\disjTau}{\mathop{\mathit{disj}_\tau}}
\newcommand{\ON}{\mathop{\mathcal{ON}}}
\newcommand{\avd}{\mathop{\mathit{\overline{d}_{av}}}}
\newcommand{\gt}{\gamma_t}

\newcommand{\cD}{\mathcal{D}}
\newcommand{\cC}{\mathcal{C}}
\newcommand{\cF}{\mathcal{F}}

\newcommand{\cT}{\mathcal{T}}

\newcounter{myFigure}
\newenvironment{myFigure}{\begin{center}}{\end{center}}
\renewcommand{\caption}[1]{\refstepcounter{myFigure}\textbf{Figure~\arabic{myFigure}. }#1}

\newcommand{\hyperedgetwo}[6]{
	\pgfmathsetmacro\Done{sqrt((#4-#1)^2+(#5-#2)^2)}
	\pgfmathsetmacro\angleone{(#2>#5)*(180+asin((#4-#1)/ 
\Done)-asin((#1-#4)/ \Done))+asin((#1-#4)/ \Done)+asin((#3-#6)/\Done)}
	
\pgfmathsetmacro\angleone{\angleone-360*(\angleone>0)-360*(\angleone>360)}
	\draw ([shift=(\angleone:#3)] #1,#2)--([shift=(\angleone:#6)]#4,#5);
	\pgfmathsetmacro\Dtwo{sqrt((#1-#4)^2+(#2-#5)^2)}
	\pgfmathsetmacro\angletwo{(#5>#2)*(180+asin((#1-#4)/ 
\Dtwo)-asin((#4-#1)/ \Dtwo))+asin((#4-#1)/ \Dtwo)+asin((#6-#3)/\Dtwo)}
	
\pgfmathsetmacro\angletwo{\angletwo-360*(\angletwo>0)-360*(\angletwo>360)}
	\draw ([shift=(\angletwo:#6)] #4,#5)--([shift=(\angletwo:#3)]#1,#2);
	\draw (#1,#2)+(\angletwo:#3)
arc(\angletwo:\angleone+360*(\angleone<\angletwo):#3);
	\draw (#4,#5)+(\angleone:#6)
arc(\angleone:\angletwo+360*(\angletwo<\angleone):#6);
}

\begin{document}

\title{Thoroughly Distributed Colorings}

\author{$^{1,2}$Wayne Goddard and \, $^2$Michael A. Henning\thanks{Research
supported in part by the South African National Research Foundation and the
University of Johannesburg}
\\ \\
$^1$School of Computing and Department of Mathematical Sciences\\
Clemson University \\
Clemson SC 29634 USA \\
\small \tt Email: goddard@clemson.edu  \\
\\
$^2$Department of Pure and Applied Mathematics \\
University of Johannesburg \\
Auckland Park, 2006 South Africa\\
\small \tt Email: mahenning@uj.ac.za  \\
}

\date{}
\maketitle

\begin{abstract}
We consider (not necessarily proper) colorings of the vertices of a graph where every color is 
thoroughly distributed, that is, appears in every open neighborhood. Equivalently,
every color is a total dominating set.
We define $\td(G)$ as the maximum
number of colors in such a coloring and $\FTD(G)$ as the fractional version thereof.
In particular, we show that every claw-free graph with minimum degree at
least~$2$ has~$\FTD(G)\ge 3/2$ and this is best possible.
For planar graphs, we show that
every triangular disc has $\FTD(G) \ge 3/2$ and this is best possible, and that 
every planar graph has $\td(G) \le 4$ and this is best possible, while we conjecture that 
every planar triangulation
has $\td(G)\ge 2$. Further,  although there are arbitrarily large
examples of connected, cubic
graphs with $\td(G)=1$, we show that
for a connected cubic graph $\FTD(G) \ge 2-o(1)$, and
conjecture that it is always at least~$2$.
We also consider the related concepts in hypergraphs.
\end{abstract}


\newpage

\section{Introduction}

We consider (not necessarily proper) colorings of the vertices of a graph where every color is 
thoroughly distributed, that is, appears in every open neighborhood. 
Chen et al.~\cite{CKTV15coupon} called this the \emph{coupon coloring problem}.

We define
$\td(G)$ as the maximum number of colors in such a coloring.
Note that a color being thoroughly distributed is equivalent to 
a color being a total dominating set, and thus the parameter is 
equivalent to the maximum number of
disjoint total dominating sets, which Cockayne et al.~\cite{CDH80total} called
the \emph{total domatic number}.
This parameter is now well studied. For example,
Zelinka~\cite{Zelinka89total} showed that there are
graphs with arbitrarily large minimum degree without two disjoint total
dominating sets. Heggernes and Telle~\cite{HeTe98} showed that the decision
problem to decide for a given graph $G$ if $\td(G) \ge 2$ is NP-complete, even
for bipartite graphs. In contrast, several researchers, such as Aram et al.~\cite{ASV12regular},
studied $\td(G)$ for a $k$-regular graph $G$; in particular, Chen et
al.~\cite{CKTV15coupon} showed that such graphs have total domatic
number at least $(1-o(1)) k/\ln k$. The related idea in hypergraphs, where every color appears in every edge, 
is called \emph{panchromatic coloring}~\cite{KoWo01panchromatic}.

In this paper we consider the fractional analogue of the parameter $\td(G)$.
We define a \emph{thoroughly distributed family} $\cF$ of a graph $G$ as a family of
(not necessarily distinct) total dominating sets of $G$.  We denote by
$r_{_\cF}$ the maximum times any vertex of~$G$ appears in~$\cF$, and
define the \emph{effective ratio} of the family $\cF$
as the ratio of the number of sets in $\cF$ to~$r_{_\cF}$.
The \emph{fractional total domatic number}
$\FTD(G)$ is then defined as the supremum of the effective ratio taken over
all thoroughly distributed families. That is,
\[
\FTD(G) = \sup_{\cF} \, \frac{ |\cF| }{ r_{_\cF}}.
\]

Like other fractional parameters, one can show that the supremum can be
achieved. For example, if we let $\cT_G$ be the hypergraph with vertex set $V(G)$
and hyperedges all total dominating sets of $G$, then $\FTD(G)$ is the
fractional matching number of $\cT_G$, and can be viewed as a linear program.
(See Chapter~1 of \cite{ScUl97fractional} for further discussion.)

We remark that there have been a few papers on the ordinary domination
equivalent. We let $\dom(G)$ denote the \emph{domatic number} of $G$, and
so $\dom(G)$ is the maximum number of disjoint dominating sets in $G$. We let
$\FD(G)$ denote the \emph{fractional domatic number}, defined analogously as
the fractional total domatic number. The fractional domatic number seems to
have been introduced by Suomela~\cite{Suomela06sleep}.
Recently, Abbas et al.~\cite{AELTWarxiv} showed that:
a $K_{1,6}$-free graph $G$ with minimum degree $\delta \ge 2$ has $\FD(G) \ge 5/2$ except for
some exceptions. See also~\cite{FYK00configurations}.

We proceed as follows.
In Section~2  we provide some preliminary results.
In Section~3  we introduce the related concepts in hypergraphs.
In Section~4 we consider claw-free graphs, and show that
every claw-free graph with minimum degree at
least~$2$ has fractional total domatic number at least~$3/2$ and this is best possible.
In Section~5 we consider planar graphs. Inter alia, we show that
every triangular disc has fractional total domatic number at
least~$3/2$ and this is best possible; that a maximal outerplanar graph
has two disjoint total dominating sets provided the order is not congruent to $2$ modulo $4$;
and in general that a planar graph has total
domatic number at most $4$ and this is best possible.
Finally in Section~6 we show that the fractional total
domatic number of a connected cubic graph is at least $2-o(1)$.
Along the way we pose several open problems.

\section{Preliminary Observations and General Properties}

Since each set in a throughly distributed family $\cF$ of a graph $G$ has size at least~$\gt(G)$, by averaging we have that $r_{_{\cF}}  \ge |\cF|
\, \gt(G) / n(G)$, implying that the effective ratio $|\cF|/r_{_{\cF}}$ is at most~$n(G)/\gt(G)$. If a 
thoroughly distributed
family $\cF$ consists of a maximum number of disjoint total dominating
sets of $G$, then $|\cF| = \td(G)$ and $r_{_{\cF}} =
1$, implying that $\FTD(G) \ge \td(G)$. We state these observations formally as
follows.

\begin{observation} \label{Ob:1}
If $G$ is an isolate-free graph of order $n$, then
\[
 \td(G)  \le \FTD(G) \le \frac{n}{\gt(G)} .
\]
\end{observation}

Equality occurs throughout the chain, for example, if $G=K_n$ and $n$ is even,
or if $G=C_n$ and $n$ is a multiple of $4$. More generally, equality
occurs in the upper bound for all complete graphs and cycles. As a gentle introduction, we determine the fractional total domatic number of a cycle. Recall that $\gt(C_n) = \lfloor n/2 \rfloor + \lceil n/4 \rceil - \lfloor n/4 \rfloor$.

\begin{observation} \label{Ob:cycle}
For $n \ge 3$, $\FTD(C_n) = n/\gt(C_n)$.
\end{observation}
\begin{Proof}
Let $G$ be the cycle $1 2 \ldots n 1$. Let $S$ be an arbitrary minimum total dominating set of $G$, and so $|S| = \gt(G)$. For $1\le i \le n$, let $S_i = S + i = \cup_{j \in S} \{i+j\}$, {where addition is taken modulo~$n$}. Each set $S_i$ is a minimum total dominating set of $G$. Let $\cF = \{S_1,S_2, \ldots, S_n\}$. Each vertex of $G$ appears in exactly~$\gt(G)$ of these sets, and therefore the effective total-ratio of the family~$\cF$ is $|\cF|/r_{_{\cF}} = n/\gt(G)$, implying that $\FTD(G) \ge n/\gt(G)$. The result now follows from the upper bound of Observation~\ref{Ob:1}.
\end{Proof}

For example, if we take $G$ as the $5$-cycle, then $\FTD(C_5) = 5/3$. However, $\td(C_5) = 1$.

\begin{observation}
\label{Ob:2}
If a graph $G$ has minimum degree $\delta \ge 1$, then $\FTD(G) \le \delta$.
\end{observation}
\begin{Proof}
Let $v$ be a vertex of minimum degree in $G$. Each set in a thoroughly distributed family $\cF$ of $G$ must contain a neighbor of the vertex~$v$. So, by the Pigeonhole Principle, at least one neighbor of~$v$ appears
in at least~$|\cF|/\delta$ sets, and so $r_{_{\cF}} \ge |\cF|/\delta$. This is true for every thoroughly distributed family $\cF$.
\end{Proof}

On the other hand, we have the following:

\begin{observation}
If a graph $G$ of order~$n$ has minimum degree $\delta \ge 1$, then $\FTD(G)
\ge n/(n-\delta+1)$.
\label{Ob:deltaLower}
\end{observation}
\begin{Proof}
Consider the collection $\cF$ of all subsets $F$ of the vertex set of exactly $n-\delta+1$
elements. Then every vertex has a neighbor in $F$; that is, $\cF$ is a thoroughly distributed family. Further, every vertex is in $\binom{n-1}{n-\delta}$ sets of $\cF$. Thus, $\cF$ has effective ratio $\binom{n}{n-\delta+1} /  \binom{n-1}{n-\delta} =
 n/(n-\delta+1)$. The result follows.
\end{Proof}

Thus, for example, by Observations~\ref{Ob:2} and \ref{Ob:deltaLower} we get the following observation.

\begin{cor}
(a) If a graph $G$ has minimum degree $\delta = 1$,  then $\FTD(G) = 1$. \\
(b) If a graph $G$ has minimum degree $\delta \ge 2$,  then $\FTD(G) > 1$.
\label{c:delta}
\end{cor}

We will show (Theorem~\ref{t:SubKn}) that there are graphs $G$ with arbitrarily
large minimum degree with $\FTD(G) < 1+ \epsilon$. Indeed, these are the
graphs that Zelinka~\cite{Zelinka89total} provided as
examples that have $\td(G) = 1$ and
arbitrarily large minimum degree.

We note that the union of two disjoint dominating sets is
a total dominating set. Thus, it is immediate that $\td(G) \ge \lfloor \dom(G)/2 \rfloor$.
But, in the case that the ordinary domatic number is odd, one can say slightly more:

\begin{theorem} \label{t:domatic}
If $G$ is an isolate-free graph, then
$\FTD(G) \ge \dom(G)/2$.
\end{theorem}
\begin{Proof}
If $D_1, \ldots, D_k$ are disjoint dominating sets, then
$\cF=  \{ \, D_i \cup D_j : 1\le i < j \le k \, \}$ is a thoroughly distributed family.
Every vertex appears in at most $k-1$ sets; so $\FTD(G) \ge \binom{k}{2}/(k-1) = k/2$.
\end{Proof}

\noindent
Equality occurs, for example, in complete graphs.

Note that $\FTD$ is monotonic, in that it cannot decrease on the addition of edges.
We next show that the fractional total domatic number of the
disjoint union behaves as expected.

\begin{theorem}
If $G$ is the disjoint union of isolate-free graphs $G_1, G_2, \ldots, G_k$,
then $\FTD(G) = \min \{ \FTD(G_1), \FTD(G_2), \ldots,  \FTD(G_k) \}$.
\label{thm:unionFTD}
\end{theorem}
\begin{Proof}
It suffices to prove this for $k=2$, as the full result follows by induction.
We use the standard notation $[s] = \{1,2,\ldots,s\}$.

For $\ell \in [2]$, let
$\cF_\ell= \{T_{\ell,1}, \ldots, T_{\ell,k_\ell}\}$ be an optimal thoroughly distributed family of $G_\ell$, 
and let $r_\ell = r_{_{\cF_\ell}}$.
Then the collection
$
\cF = \{ \, T_{1,i} \cup T_{2,j} : i \in [k_1], \,j \in [k_2] \,\}
$
is a thoroughly distributed family of~$G$.
If a vertex $v$ of $G_1$ appears in $r_v$ sets
of $\cF_1$, then $v$ appears in $r_v \times k_2$ sets of~$\cF$,
and similarly with vertices of $G_2$.
Thus, $r_{_{\cF}} = \max \{ r_1k_2, r_2k_1 \} $.
So,
\[
\FTD(G) \ge  \frac{k_1k_2}{\max\{ r_1k_2, r_2k_1\} } =
\min \left\{ \frac{k_1}{r_1},  \frac{k_2}{r_2} \right\} = \min \{ \FTD(G_1), \FTD(G_2) \} .
\]

Conversely, let $\cF^*= \{ U_{1}, \ldots, U_{k} \} $ be an optimal thoroughly distributed family of $G$.
Then, for $\ell \in [2]$, the collection
$
 \cF^*_\ell = \{\, U_i \cap V(G_\ell) : i \in [k] \,\}
$
is a thoroughly distributed family of $G_\ell$. Further,
if each vertex appears at most $r$ times in $\cF^*$, then
$r_{_{\cF^*_{\ell}}} \le r$. Thus
\[
  \FTD(G_\ell) \ge \frac{k}{r} = \FTD(G).
\]
The two inequalities combined give the result.
\end{Proof}

\section{Hypergraphs}

We observed earlier that the fractional total domatic number of $G$ is also
the fractional matching number of the hypergraph $\cT_G$. But hypergraphs
also provide a more general setting for the parameter.

Recall that a subset $T$ of vertices in a hypergraph $H$ is a \emph{transversal} (also
called \emph{vertex cover},  \emph{hitting set} or \emph{blocking set})
if $T$ has a nonempty intersection with every edge of $H$. The
\emph{transversal number} $\tau(H)$ of $H$ is the minimum size of a transversal
in $H$. See for example~\cite{Berge79survey,BuHeTu12,BuHeTuYe14}.

A hypergraph $H$ is $2$-\emph{colorable} if there is a $2$-coloring of the
vertices such that each hyperedge contains two vertices of distinct colors. In
other words, there is no monochromatic hyperedge. So, the question of
when a hypergraph has two disjoint transversals is the same as
whether the hypergraph has a $2$-coloring (also known as Property~B), or whether
a design has a blocking set. More generally, Kostochka and Woodall~\cite{KoWo01panchromatic}
defined a panchromatic $k$-coloring of a hypergraph as a coloring with $k$ colors
such that every hyperedge contains each color. This is equivalent to a partition
into $k$ disjoint transversals. We denote by
$\disjTau(H)$ the \emph{disjoint transversal number} of a
hypergraph $H$, which is the maximum number of disjoint
transversals in $H$.

Analogous to the fractional total domatic number, one can define the fractional
disjoint transversal number.
A \emph{transversal family} $\cF$ of a hypergraph $H$ is a family of
transversals of~$H$. Given a hypergraph~$H$ and a transversal family
$\cF$, we define the \emph{effective transversal-ratio} of the
family~$\cF$ as the ratio of the number of sets in $\cF$ over the maximum
times~$r_{_\cF}$ any element appears in~$\cF$. The \emph{fractional disjoint
transversal number} $\FdT(H)$ is the supremum of the effective
transversal-ratio taken over all transversal
families. That is,
\[
\FdT(H) = \sup_{\cF} \, \frac{ |\cF| }{\, r_{_{\cF}}} .
\]

Analogous to Observation~\ref{Ob:1}, we have the following bounds on the
fractional disjoint transversal number.

\begin{observation}  \label{ob:hypergraph}
For every isolate-free hypergraph $H$ of order $n$,
\[
 \disjTau (H) \le \FdT(H) \le \frac{n}{\tau(H)} .
\]
\end{observation}

For example, if $H$ is the complete $k$-uniform hypergraph of order $n$,
then $\tau(H) = n-k+1$, and $\FdT(H) = n/(n-k+1)$, by considering
the collection of all $n-k+1$ element subsets as transversal family.

As another example, if we take $H$ to be the Fano plane and let $\cF$ be
the transversal family consisting of the seven edges of $H$, then
$|\cF| = 7$ and $r_{_{\cF}} = 3$, implying that $\FdT(H) \ge 7/3$, with equality by
the above observation. However, $\disjTau (F_7) = 1$.

And similar to Observation~\ref{Ob:deltaLower} we have:

\begin{observation}
If every vertex of hypergraph $H$ is in an edge of size at least~$a$, then $\FdT(H) \ge n/(n-a+1)$.
\label{ob:hyperLower}
\end{observation}
\begin{Proof}
{Consider the collection $\cF$ of all subsets $F$ of the vertex set of exactly $n-a+1$ elements. Then every edge contains a vertex in $F$; that is, $\cF$ is a transversal family. Further, every vertex is in $\binom{n-1}{n-a}$ sets of $\cF$. Thus, $\cF$ has effective ratio $\binom{n}{n-a+1} /  \binom{n-1}{n-a} = n/(n-a+1)$. The result follows.} \end{Proof}

Analogous to Theorem~\ref{thm:unionFTD}, there is the following result on the
fractional disjoint transversal number of the disjoint union of hypergraphs.

\begin{theorem}
If $H$ is the disjoint union of isolate-free hypergraphs $H_1, H_2, \ldots, H_k$,
then $\FdT(H) = \min \{ \FdT(H_1), \FdT(H_2), \ldots, \FdT(H_k) \}$.
\label{thm:unionHyper}
\end{theorem}

\subsection{Connection}

Associated with a graph $G$, one can define the \emph{open neighborhood
hypergraph}, abbreviated
ONH, of $G$ as the hypergraph $\ON(G)$ whose vertex set is $V(G)$ and whose
hyperedges are the open neighborhoods of vertices in $G$. Thus, if $H =
\ON(G)$, then $V(H) = V(G)$ and $E(H) = \{ \, N_G(x) : x \in V(G) \,\}$.

As Thomass\'{e} and Yeo~\cite{ThYe07} observed, a
total dominating set in $G$ is a transversal in $\ON(G)$, and conversely.
More generally,
a thoroughly distributed family of a graph $G$ is a transversal family of
$\ON(G)$ and conversely. Thus, the fractional total domatic number of
an isolate-free graph is precisely the fractional disjoint transversal number
of its ONH.

\begin{observation}
For every isolate-free graph $G$, $\FTD(G) = \FdT(\ON(G))$.
\label{ob:connection}
\end{observation}

One can also construct the incidence graph $I(H)$ of a hypergraph $H$. For example,
the incidence graph of the complete $r$-uniform hypergraph $H$ on $n$ vertices with
$n\ge 2r-1$ is the standard example of a graph that does not have
two disjoint total  dominating sets but has arbitrarily large minimum degree,
as shown by Zelinka~\cite{Zelinka89total}.

In general, every connected bipartite graph has an ONH that consists of two components
(see \cite{HeYe08}).
For example, consider the Heawood graph $G_{14}$.
The ONH of the Heawood graph consists of two disjoint copies of the Fano plane, $F_7$. See
Figure~\ref{f:Heawood}. Conversely,
the Heawood graph is the incidence graph of the Fano
plane.
We note that $\gt(G_{14}) = 6$ and
$\tau(\ON(G_{14})) = 2\tau(F_7) = 6$.
Further, by
Observation~\ref{ob:connection} and Theorem~\ref{thm:unionHyper}, we have $\FTD(G_{14}) =
\FdT(\ON(G_{14})) = \FdT(F_7 \cup F_7) = \FdT(F_7) = 7/3$. In contrast,
the Heawood
graph does not have two disjoint total dominating sets; that is, $\td(G_{14}) = 1$.

\begin{myFigure}
\begin{tikzpicture}

 \begin{scope}[xshift=0cm,yshift=0cm,thick,scale=.65]
 \tikzstyle{every node}=[circle, draw, fill=black!0, inner sep=0pt,minimum
width=.16cm]
 \draw(0,0) { 
    +(2.09,4.29) -- +(3.02,4.07)
    +(3.02,4.07) -- +(3.76,3.48)
    +(3.76,3.48) -- +(4.18,2.62)
    +(4.18,2.62) -- +(4.18,1.67)
    +(4.18,1.67) -- +(3.76,0.81)
    +(3.76,0.81) -- +(3.02,0.21)
    +(3.02,0.21) -- +(2.09,0.00)
    +(2.09,0.00) -- +(1.16,0.21)
    +(1.16,0.21) -- +(0.41,0.81)
    +(0.41,0.81) -- +(0.00,1.67)
    +(0.00,1.67) -- +(0.00,2.62)
    +(0.00,2.62) -- +(0.41,3.48)
    +(0.41,3.48) -- +(1.16,4.07)
    +(1.16,4.07) -- +(2.09,4.29)
    +(3.02,4.07) -- +(0.00,1.67)
    +(2.09,4.29) -- +(3.76,0.81)
    +(3.76,3.48) -- +(2.09,0.00)
    +(4.18,2.62) -- +(0.41,3.48)
    +(4.18,1.67) -- +(0.41,0.81)
    +(3.02,0.21) -- +(0.00,2.62)
    +(1.16,0.21) -- +(1.16,4.07)
    +(2.09,4.29) node[fill=black!100]{}  
    +(1.16,4.07) node[fill=black!100]{}  
    +(0.41,3.48) node[fill=black!100]{}
    +(0.00,2.62) node[fill=black!100]{}
    +(0.00,1.67) node[fill=black!100]{}
    +(0.41,0.81) node[fill=black!100]{}
    +(1.16,0.21) node[fill=black!100]{}
    +(2.09,0.00) node[fill=black!100]{}
    +(3.02,0.21) node[fill=black!100]{}
    +(3.76,0.81) node[fill=black!100]{}
    +(4.18,1.67) node[fill=black!100]{}
    +(4.18,2.62) node[fill=black!100]{}
    +(3.76,3.48) node[fill=black!100]{}
    +(3.02,4.07) node[fill=black!100]{}
    +(2.1,4.7) node[rectangle, draw=white!0, fill=white!100]{\small $a$}
    +(3.35,4.35) node[rectangle, draw=white!0, fill=white!100]{\small $1$}
    +(4.15,3.65) node[rectangle, draw=white!0, fill=white!100]{\small $b$}
    +(4.55,2.7) node[rectangle, draw=white!0, fill=white!100]{\small $2$}
    +(4.6,1.7) node[rectangle, draw=white!0, fill=white!100]{\small $c$}
    +(4.15,0.75) node[rectangle, draw=white!0, fill=white!100]{\small $3$}
    +(3.275,-0.1) node[rectangle, draw=white!0, fill=white!100]{\small $d$}
    +(2.1,-0.45) node[rectangle, draw=white!0, fill=white!100]{\small $4$}
    +(1.05,-0.2) node[rectangle, draw=white!0, fill=white!100]{\small $e$}
    +(1.05,4.5) node[rectangle, draw=white!0, fill=white!100]{\small $7$}
    +(0.05,0.75) node[rectangle, draw=white!0, fill=white!100]{\small $5$}
    +(0.05,3.65) node[rectangle, draw=white!0, fill=white!100]{\small $g$}
    +(-0.375,1.7) node[rectangle, draw=white!0, fill=white!100]{\small $f$}
    +(-0.375,2.7) node[rectangle, draw=white!0, fill=white!100]{\small $6$}
    +(2.1,-1.85) node[rectangle, draw=white!0, fill=white!100]{(a) $G_{14}$}    
  };
  \end{scope}

\begin{scope}[xshift=5cm,yshift=0cm,scale=0.4]
\fill (0,0) circle (0.2cm); \fill (4,0) circle (0.2cm); \fill (8,0)
circle (0.2cm); \fill (2,3.464) circle (0.2cm); \fill (6,3.464)
circle (0.2cm); \fill (4,6.928) circle (0.2cm); \fill (4,2.3093)
circle (0.2cm); \hyperedgetwo{0}{0}{0.5}{8}{0}{0.6};
\hyperedgetwo{4}{6.928}{0.5}{0}{0}{0.6};
\hyperedgetwo{4}{6.928}{0.6}{8}{0}{0.5};
\hyperedgetwo{0}{0}{0.4}{6}{3.464}{0.45};
\hyperedgetwo{4}{6.928}{0.4}{4}{0}{0.45};
\hyperedgetwo{8}{0}{0.4}{2}{3.464}{0.45}; \draw (4,2.3093) circle
(1.9523cm); \draw (4,2.3093) circle (2.6523cm);
\node at (4,8) {\small $1$};
\node at (7.05,4) {\small $2$};
\node at (1,4) {\small $3$};
\node at (0,-1.2) {\small $7$};
\node at (4,-1.2) {\small $5$};
\node at (8,-1.2) {\small $4$};
\node at (4,3.35) {\small $6$};
\node at (9,-3) {(b) $\ON(G_{14})$};
\end{scope}

\begin{scope}[xshift=9cm,yshift=0cm,scale=0.4]
\fill (0,0) circle (0.2cm); \fill (4,0) circle (0.2cm); \fill (8,0)
circle (0.2cm); \fill (2,3.464) circle (0.2cm); \fill (6,3.464)
circle (0.2cm); \fill (4,6.928) circle (0.2cm); \fill (4,2.3093)
circle (0.2cm); \hyperedgetwo{0}{0}{0.5}{8}{0}{0.6};
\hyperedgetwo{4}{6.928}{0.5}{0}{0}{0.6};
\hyperedgetwo{4}{6.928}{0.6}{8}{0}{0.5};
\hyperedgetwo{0}{0}{0.4}{6}{3.464}{0.45};
\hyperedgetwo{4}{6.928}{0.4}{4}{0}{0.45};
\hyperedgetwo{8}{0}{0.4}{2}{3.464}{0.45}; \draw (4,2.3093) circle
(1.9523cm); \draw (4,2.3093) circle (2.6523cm);
\node at (4,8) {\small $a$};
\node at (7,4) {\small $b$};
\node at (1.1,4) {\small $g$};
\node at (0,-1.2) {\small $e$};
\node at (4,-1.2) {\small $c$};
\node at (8,-1.2) {\small $f$};
\node at (4,3.5) {\small $d$};
\end{scope}

\end{tikzpicture}
\\
\caption{The Heawood graph and its open neighborhood hypergraph}
\label{f:Heawood}
\end{myFigure}

The \emph{subdivision graph} of a graph $G$, denoted $S(G)$, is the graph
obtained from $G$ by subdividing every edge of $G$ exactly once.

\begin{theorem}
For $n \ge 3$, $\FTD( S(K_n) ) =  n/(n-1)$.
\label{t:SubKn}
\end{theorem}
\begin{Proof}
Let $G = S(K_n)$ with vertex set $V(K_n) \cup E(K_n)$. Then, the ONH of $G$ consists of two components. One, say $H_1$, has vertex set $V(K_n)$ and one, say $H_2$, has vertex set $E(K_n)$.
The hypergraph $H_1$ is $2$-uniform and is the graph $K_n$,
while hypergraph $H_2$ is $(n-1)$-uniform, $2$-regular with $n$ edges.
By Observation~\ref{ob:hyperLower}, $\FdT(H_i) \ge n/(n-1)$ for $i \in [2]$. A transversal of $
H_1 \cong K_n$ is a vertex cover in $H_1$, implying that $\tau(H_1) = n-1$. Thus, by
Observation~\ref{ob:hypergraph}, $\FdT(H_1) \le n(H_1)/\tau(H_1) = n/(n-1)$.
Consequently, $\FdT(H_1) = n/(n-1)$.
The result follows by Theorem~\ref{thm:unionHyper}.
\end{Proof}

Note that $S(K_n)$ is also the incidence graph of $K_n$, if $K_n$
is thought of as a $2$-uniform hypergraph. More generally, we have the
following result.

\begin{theorem}
Let $G_{k,n} = I(H )$ be the incidence graph of the complete $k$-uniform
hypergraph. Then, $\FTD( G_{k,n} ) =  n/(n-k+1)$ for $n\ge 2k-1$.
\end{theorem}

\section{Claw-Free Graphs}

We show next that a claw-free graph $G$ with minimum degree at least~$2$ has
$\FTD(G) \ge 3/2$. We will need the following results.

Recall that if $G$ is a graph, $S$ a subset of vertices of $V(G)$, and $v$ a
vertex of $G$, then an $S$-\emph{private neighbor of $v$} is a neighbor of $v$
that is not a neighbor of any other vertex of~$S$. Cockayne, Dawes, and
Hedetniemi~\cite{CDH80total} established the following property.

\begin{proposition}{\rm (\cite{CDH80total})}
\label{p:properties} If
$S$ is a minimal total dominating set in graph $G$, then every vertex in $S$
has an $S$-private neighbor.
\end{proposition}

We will also need the result that every graph has two disjoint sets, one
total dominating and one ``half-total'' dominating:

\begin{theorem}{\rm (\cite{HS08disjoint})}
If $G$ is a graph with minimum degree at least~$2$, then, except for the $5$-cycle,
$G$ has disjoint dominating and total dominating sets.
\label{t:domTdomP}
\end{theorem}

We can now establish a lower bound on the fractional total
domatic number of a claw-free graph with minimum degree at least~$2$.

\begin{theorem}
If $G$ is a claw-free graph with $\delta \ge 2$, then $\FTD(G) \ge 3/2$.
\label{t:claw2}
\end{theorem}
\begin{Proof}
As observed earlier, $\FTD(C_5) = 5/3$.  Hence, we may assume
that $G$ is not a $5$-cycle. By
Theorem~\ref{t:domTdomP}, the graph $G$ therefore has a partition $R$ and $B$
such that $R$ is a total dominating set and $B$~is a dominating set.
By transferring vertices from $R$ to $B$ if necessary, we
can assume $R$ is a minimal total dominating set.

We will construct three total dominating sets: one is $R$,
and the other two $B_1$ and $B_2$ are supersets of $B$.
We need to partition $R$ between $B_1$ and $B_2$. If this can be done, then the
resulting thoroughly distributed family $\cF_{_G} = \{R,B_1,B_2\}$ has the property
that every vertex of $G$ appears in exactly two sets of $\cF_{_G}$, so
that the effective ratio of $\cF_{_G}$ is~$3/2$ and hence $\FTD(G) \ge
3/2$.

Let $B'$ be the isolated vertices of $B$. The only vertices that do not
have a neighbor in~$B$ are $B'$.
Now, for every vertex $x$ in $B'$ choose two neighbors $x_1$
and $x_2$ (necessarily
in~$R$) subject to the constraint that if possible one is an $R$-private
neighbor of the other.

Now, construct an auxiliary graph $A$ with vertex set $R$ and
two vertices in $A$ are adjacent if and only if they are the $\{x_1,x_2\}$
pair of some $x\in B'$. Observe that if $A$ is bipartite
then we are done---that is, the $2$-coloring is the desired partition
of~$R$---since every vertex in $B'$ is adjacent to
a vertex in both partite sets.

So assume $A$ is not bipartite. Let $C \colon v_1v_2 \ldots v_k v_1$ be a
shortest odd cycle in $A$. For $i \in [k]$, let $u_i$ be the vertex of $B'$
associated with the pair~$\{v_i,v_{i+1}\}$ (where addition is taken
modulo~$n$). Thus, in $G$, $v_1u_1v_2u_2 \ldots v_ku_kv_1$ is a cycle. By
claw-free-ness, every vertex in~$R$ has at most two neighbors in $B'$. In
particular, for each vertex $v_i$ of $C$, its neighbors in~$B'$ are the two
vertices $u_{i-1}$ and $u_{i}$ (where addition is taken modulo~$n$). This
implies, by claw-free-ness, that the $R$-private neighbors of a vertex $v \in V(C)$ all belong to
$R$.

Because $C$ is odd, and $R$-private neighbors are unique or mutual, there must
be a vertex $v$ in $C$ whose $R$-private neighbor $w$ is not on $C$. Renaming
vertices, if necessary, we may assume that $v_1$ is such a vertex in $C$. We
now consider the neighbors of $v_1$ in $B'$, namely $u_1$ and $u_k$. By the
claw-free-ness of $G$, $u_1w$ or $u_kw$ is an edge of $G$. We may assume that
$u_1w$ is an edge. Then by the choice of $x$'s neighbors, it must be that $v_2$
is another $R$-private
neighbor for $v_1$ but it is on $C$, a contradiction.
\end{Proof}

We note that the lower bound of Theorem~\ref{t:claw2} is somewhat best possible:
the graphs $K_3$ and
$C_6$ have fractional total domatic number exactly $3/2$. Further, there are graphs
with claws, such as the subdivision $S(K_4)$, that have smaller fractional total
domatic number (see Theorem~\ref{t:SubKn}).
Further, there are arbitrarily large connected $K_{1,4}$-free graphs with fractional total domatic number exactly $3/2$. For example, take {an even number $k \ge 2$ of} disjoint copies of the $6$-cycle, {let $(R,B)$ denote a bipartition of the resulting graph $kC_6$,} and add a matching between the vertices of $R$ to make the graph connected. Since $C_6$ has three disjoint dominating sets, it follows from Theorem~\ref{t:domatic} that the resultant graph $G$ has $\FTD(G)\ge 3/2$; on the other hand, every total dominating set {of $G$} must contain at least {two vertices of $R$ from each original $6$-cycle in order to totally dominate the three vertices of $B$ that belong to that cycle, implying that $\gt(G) \ge 2|R|/3$. Hence, if $\cF$ is a thoroughly distributed family of $G$, then, by averaging, we have that $r_{_{\cF}}  \ge |\cF| 
\, \gt(G) / |R|$, implying that the effective ratio $|\cF|/r_{_{\cF}}$ is at most~$|R|/\gt(G) \le 3/2$ and hence $\FTD(G) \le 3/2$. Consequently, $\FTD(G) = 3/2$.}

However, we believe that the lower bound of Theorem~\ref{t:claw2} should be improvable asymptotically.
Perhaps it is true that if $G$ is a connected, claw-free graph with $\delta \ge 2$,
then $\FTD(G) \ge 2 - o(1)$ and/or there is a partition $(T_1,T_2)$ of the vertex set
such that every vertex except possibly two has a neighbor in both $T_1$ and $T_2$.
And maybe, if we require that $\delta(G) \ge 3$, then two disjoint
total dominating sets are guaranteed in $G$.

\section{Planar Graphs}

In this section, we consider the fractional total domatic number of planar
graphs. Of course in general there are no lower bounds. So we focus on
dense graphs.

\subsection{Triangulated Discs}

A \emph{triangulated disc} is a (simple) planar graph all of whose faces are triangles,
except possibly for the outer face. Matheson and Tarjan~\cite{MaTa96}
showed that if $G$ is a triangulated disc, then $\dom(G) \ge 3$. Hence, as an
immediate consequence of Theorem~\ref{t:domatic} and the Matheson-Tarjan
result, we have the following lower bound.

\begin{theorem}
\label{t:disc}
If $G$ is a triangulated disc, then $\FTD(G) \ge \frac{3}{2}$.
\end{theorem}
We remark that the lower bound of Theorem~\ref{t:disc} is tight, in
that there exist triangulated discs $G$ of arbitrarily large order satisfying
$\FTD(G) = 3/2$. For example, consider the triangulated disc $G$
illustrated in Figure~\ref{f:tDisc}, where the shaded area consists of any
maximal planar graph (or, equivalently, triangulation). Let $S$ be the set of
three vertices on the outer face of $G$ that have degree at least~$4$.
Since each set in a thoroughly distributed family~$\cF$ contains at least two
vertices of $S$, by averaging there is a vertex in $S$ that belongs to at
least~$2|\cF|/3$ sets in $\cF$, and so the effective ratio $|\cF|/r_{_{\cF}}$ is at most~$3/2$.
Consequently, by Theorem~\ref{t:disc}, $\FTD(G) = 3/2$.

\begin{myFigure}
\includegraphics{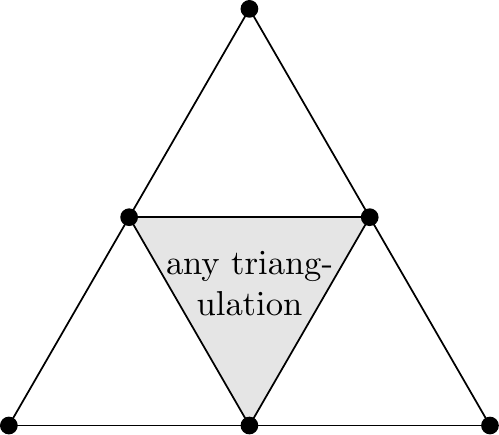}
\\
\vspace{0.2cm}
\caption{A triangulated disc $G$ with $\FTD(G) = 3/2$}
\label{f:tDisc}
\end{myFigure}

There are two extremal examples of triangulated discs: maximal planar graphs (where
the outer face is a triangle) and maximal outerplanar graphs (where the
outer face contains all vertices). We consider these next.

\subsection{Maximal Outerplanar Graphs}

The total domination number of maximal outerplanar graphs has
been studied by several authors. In particular,
Dorfling et al.~\cite{DHJ16outerplanar} showed that, except for two
exceptions, every maximal outerplanar graph with order $n\ge 5$ has total
domination number at most~$2n/5$.
Since a maximal outerplanar graph has minimum degree~$2$, the most we
can hope for is two disjoint total dominating sets (Observation~\ref{Ob:2}). In this subsection we
investigate when this occurs.

We will use the following construction.
For a maximal outerplanar graph $G$ of order $k\ge 3$, define the
graph $M(G)$ as follows. Start with $G$ and,
for every edge $e=uv$ on the outer boundary, add a new vertex $w_e$ with
neighbors $u$ and $v$. Note that $M(G)$ has order~$2k$ and is maximal
outerplanar. For example, $M(K_3)$ is the Haj\'{o}s graph or $3$-sun.
Another example of such a graph $M(G)$ is shown in Figure~\ref{f:graphMG}, where
$G$ is the maximal outerplanar graph induced by the darkened vertices.

\begin{myFigure}
\includegraphics{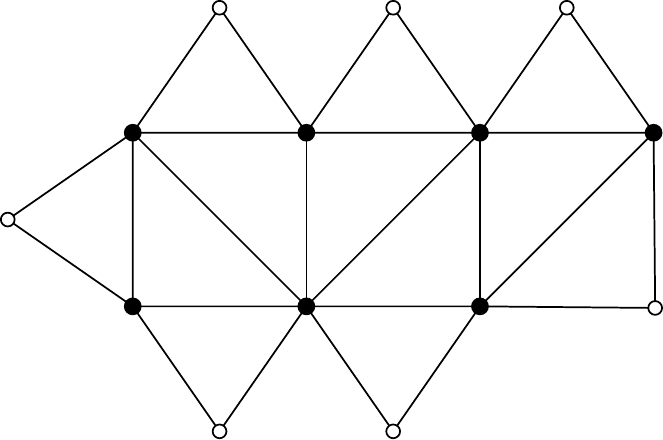}
\\
\vspace{0.2cm}
\caption{A graph $M(G)$}
\label{f:graphMG}
\end{myFigure}

\begin{observation}
If $G$ is a maximal outerplanar graph of \textbf{odd} order~$k \ge 3$,
then $M(G)$ does not have two disjoint total dominating sets.
\end{observation}
\begin{Proof}
The subhypergraph of $\ON( M(G) )$ induced by the open neighborhoods
of all the~$w_e$ is a $2$-uniform hypergraph isomorphic to an $k$-cycle.
Since $k$ is odd, such a cycle is not $2$-colorable, which means that
$\ON( M(G) )$ is not $2$-colorable. That is, $M(G)$ does not have
two disjoint total dominating sets.
\end{Proof}

For the proof of the next result we will need to recall a few concepts.
The \emph{weak dual} of a triangulated disc is the graph that has
a vertex for every bounded face of the embedding, and an edge for every pair of
adjacent bounded faces. It is known that if $G$ is a maximal
outerplanar graph of order $n$, then the weak dual graph of $G$ is a tree of order $n-2$ and maximum degree at most~$3$.
Further, recall that there is a canonical bijection between the edges of a
planar graph and the edges of its dual.

\begin{theorem}
If $G$ is a maximal outerplanar graph of order $n\ge 4$ with $n$ not congruent
to~$2$
modulo~$4$, then $G$ has two disjoint total dominating sets.
\label{t:outerp}
\end{theorem}
\begin{Proof}
Consider an embedding of the maximal outerplanar graph $G$ in the plane. Let
$C$ denote the outer hamiltonian cycle of $G$. The result is
trivial if $n$ is a multiple of $4$, since the cycle $C$ has the property.

If $n$ is congruent to $1$ modulo $4$, then take the cycle $C$ and
consider a vertex~$v$ that has a neighbor $w$ on $C$ of degree~$2$ in $G$.
Give $v$ and both its neighbors on $C$ the same color; then alternate colors in
pairs along the cycle $C$. The only vertex whose neighbors on the cycle are the
same color is $v$; but $v$ is adjacent to $w$ and $w$'s other neighbor,
which receive different colors. Thus, every vertex sees both colors, implying
that $G$ has two disjoint total dominating sets.

So assume $n$ is congruent to $3$ modulo $4$, and so $n \ge 7$. Now, we claim
that a maximal outerplanar graph has a chord
$e=uv$ such that the vertices $u$ and $v$ are distance $3$ or $4$ on the
cycle~$C$. For, consider
the weak dual $D$ of $G$. Recall, $D$ is a tree of order at least~$5$ and maximum degree at most~$3$.
Let $t$ be a vertex of $D$ that is not a leaf but has exactly one non-leaf neighbor $t'$. Then the desired chord is the edge $uv$ that is the dual of the edge~$tt'$, noting that if $t$ has degree~$2$ in $D$, then $u$ and $v$ are
distance $3$ on the cycle~$C$, while if $t$ cannot be chosen to have degree~$2$, and therefore $t$ has degree~$3$ in $D$, then $u$ and $v$ are distance~$4$ on~$C$. For the coloring of~$G$, give $u$ red, then the
next two on the cycle blue, the next two red, and so on. Note that $v$ gets
color red, and so every vertex sees both colors.
\end{Proof}

The above result is best possible, in that for each $n$ congruent
to~$2$ modulo $4$ there is a maximal outerplanar graph without two disjoint total dominating sets, namely the graphs $M(G)$ defined above.

\begin{theorem}
If $G$ is a maximal outerplanar graph of order $n\ge 6$ congruent to $2$
modulo~$4$, then $\FTD(G) \ge 2n/(n+2)$, and this is best possible.
\label{t:outerp2}
\end{theorem}
\begin{Proof}
Since $G$ is hamiltonian, $\FTD(G) \ge \FTD(C_n) = 2n/(n+2)$, by Observation~\ref{Ob:cycle}.

To show that this bound is best possible, consider a maximal planar
graph $F$ of order $k=n/2$ and the graph $M(F)$ defined
above. Since $n\ge 6$ is congruent to $2$ modulo~$4$, we note that $k \ge 3$ is odd.
Consider a total dominating set $T$ of $M(F)$. Then, as above,
the set $T$ must contain one of every consecutive pair on the outer cycle of $F$,
so that $T$ contains at least $(k+1)/2$ vertices of $F$.
Thus if $\cF$ is a thoroughly distributed family of~$M(F)$, by averaging
at least one vertex in $V(F)$ appears in at least
$|\cF| \left( \frac{k+1}{2} \right) / |V(F)|$ sets, and so $r_{_{\cF}}
\ge \left( \frac{k+1}{2k} \right)  |\cF|$. Therefore,
\[
\frac{ |\cF| }{\, r_{_{\cF}}} \le \frac{2k}{k+1} = \frac{2n}{n+2}.
\]
As we have the matching lower bound, it follows that $\FTD( M(F) ) = 2n/(n+2)$.
\end{Proof}

We remark that the graphs $M(F)$ are not the only maximal outerplanar
graphs $G$ that have $\td (G)=1$ and $\FTD(G)=2n/(n+2)$. Nevertheless,
it can be shown that all examples share the feature of having an independent set
of size $n/2$.

\subsection{Upper Bounds for Triangulations}

In this section we consider (simple) planar triangulations or equivalently
maximal planar graphs. We denote the average degree of a graph $G$
by $\avd(G)$.

\begin{lemma}
\label{l:triangulation}
If $G$ is a (not necessarily planar) triangulation of order at least~$4$, then $\FTD(G) \le
\avd(G) - 1$.
\end{lemma}
\begin{Proof}
Let $G$ be a triangulation with vertex set $V$ of order~$n \ge 4$.
Consider a thoroughly distributed family $\cF = \{F_1, F_2, \ldots, F_k\}$ of
$G$. For every $i \in [k]$ and
every pair $u$, $v$ of vertices of $G$, define a weight function $g_i(u,v)$
as follows:\\
\textbullet{}
$g_i(u,v) = 0$ if $u \notin F_i$, or if $v$ is not a neighbor of $u$;\\
\textbullet{}
$g_i(u,v) = 1/2$ if $u$ and $v$ lie in a triangle with some vertex $w$ where $u,w \in F_i$; and \\
\textbullet{}
$g_i(u,v) = 1$ if $u \in F_i$, $v$ is a neighbor of $u$, but $u$ and $v$ have no
common neighbor in $F_i$.

Consider any specific $i \in [k]$ and specific vertex $u \in F_i$. Since $F_i$
is a total dominating set of $G$, the vertex $u$ is guaranteed to have a
neighbor $w$ in $F_i$. Thus, since $G$ has minimum degree at least $3$,
there are least two choices for $v$ such that $g_i(u,v)=1/2$. The
remaining neighbors $v$ of $u$ all satisfy $g_i(u,v) \le 1$.
It follows that $\sum_{v \in V} g_i(u,v) \le \deg u-1$. More generally, since
$u$ appears in at most~$r_{_\cF}$ sets $F_i$,
\begin{equation}
\label{Eq1}
\sum_{u \in V} \left( \sum_{i=1}^k \sum_{v \in V} g_i(u,v) \right) \le  \sum_{u
\in V} r_{_\cF} (\deg u-1).
\end{equation}

On the other hand, consider any specific vertex $v \in V$. For each $i \in
[k]$, the vertex $v$ has a neighbor $u$ in $F_i$. If $u$ and $v$ are in a
triangle with some other vertex $w$ in $F_i$, then
$g_i(u,v) = g_i(w,v) = 1/2$; otherwise, $g_i(u,v) = 1$. Thus, $\sum_{u \in V}
g_i(u,v) \ge 1$ for each $i \in [k]$, and so
\begin{equation}
\label{Eq2}
\sum_{v \in V} \left(\sum_{i=1}^k \sum_{u \in V} g_i(u,v) \right) \ge \sum_{v
\in V} \left(\sum_{i=1}^k 1 \right) =  n \cdot k.
\end{equation}

By Inequalities~(\ref{Eq1}) and~(\ref{Eq2}), it follows that
\[
k \le r_{_\cF} \cdot  \frac{1}{n} \sum_{u \in V} (\deg u-1)  = r_{_\cF} (\, \avd(G) - 1) .
\]
Rearranged, this says that $|\cF| / r_{_{\cF}}  = k / r_{_\cF} \le \avd(G) - 1$, as required.
\end{Proof}

If $G$ is a planar triangulation of order~$n$, then $\avd(G) = 6 -
12/n$. Thus, as an immediate consequence of Lemma~\ref{l:triangulation}, we
have the following upper bounds.

\begin{theorem}
\label{c:triangulation}
A planar graph has total domatic number at most~$4$ and fractional
total domatic number at most $5 - \frac{12}{n}$.
\end{theorem}

There do exist planar graphs $G$ with $\td(G) = 4$.  Computer
search shows that the smallest such graph has order $16$.
For example, take the truncated tetrahedron and add
a vertex inside each hexagonal face that is joined to all vertices on the boundary.
Illustrated below (see Figure~\ref{f:td4}) is a spanning subgraph
of this that still has four disjoint total dominating sets: the vertices labelled $i$
form a total dominating set for each $i \in [4]$.

\begin{myFigure}
\includegraphics{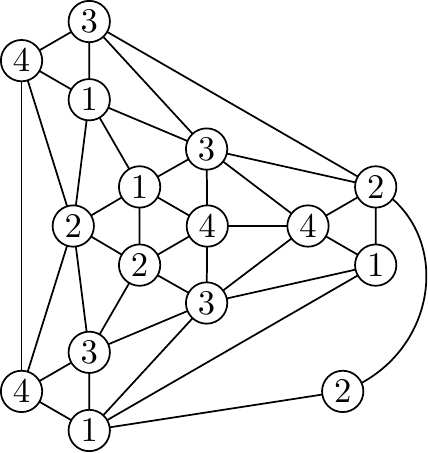}
\\
\vspace{0.2cm}
\caption{A planar graph $G$ with $\td(G) = 4$}
\label{f:td4}
\end{myFigure}

For a general construction, one can take multiple copies of the above graph
(or the graph shown in Figure~\ref{f:bigFTD})
and connect up arbitrarily. (We remark that while this construction only
produces graphs with restricted orders, it is possible to achieve intermediate
orders.)

Note that Lemma~\ref{l:triangulation} applies on all surfaces. For
example, the average degree of a maximal toroidal graph is~$6$.
It follows that:

\begin{theorem}
A toroidal graph has total domatic number at most~$5$.
\end{theorem}

We remark that there do exist toroidal graphs $G$ with $\td(G) = 5$. Such an
example is illustrated in Figure~\ref{f:torus}, where the vertices labelled $i$ form a
total dominating set of $G$ for each $i \in [5]$.

\begin{myFigure}
\includegraphics{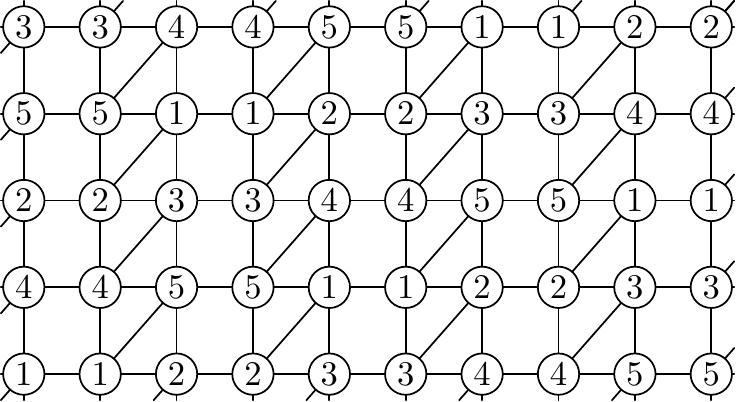}
\\
\vskip 0.2cm
\caption{A toroidal graph $G$ with $\td(G) = 5$}
\label{f:torus}
\end{myFigure}

By Theorem~\ref{c:triangulation}, a planar graph of order~$n$ has fractional
total domatic number at most $5 - \frac{12}{n}$.
The smallest planar graph $G$ for which $\FTD(G) > 4$ is illustrated in Figure~\ref{f:bigFTD}.
\begin{myFigure}
\includegraphics{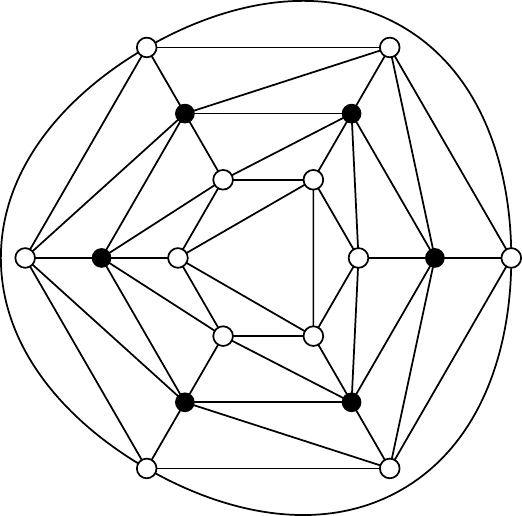}
\\
\vskip 0.2cm
\caption{A planar graph $G$ with $\FTD(G) = 21/5$}
\label{f:bigFTD}
\end{myFigure}

\subsection{A Construction}

As examples for lower bounds, we
will need the following constructions.

For a graph $G$ $2$-cell embedded on a surface, the graph
$T(G)$ is obtained from the graph~$G$
by adding in each face of $G$ a new vertex adjacent to every vertex on the
face. In the case that $G$ is planar, this is sometimes called the
\emph{Kleetope} of $G$. We will call $G$ the \emph{base} graph.

In order to prove our next result, we will need the following structural lemma
of triangulations.

\begin{lemma}
\label{lem:match}
If $G$ is a (not necessarily planar) triangulation of order at least~$4$ with set of vertices $V$
and set of faces $F$, then there is an injection from $V$ to $F$ such that
 every vertex is matched with a face on which it lies.
\end{lemma}
\begin{Proof}
We use Hall's theorem. Define for
a set $S$ of vertices, the number $f(S)$ of faces incident with at least
one vertex of $S$. It suffices to show that
$f(S)\ge |S|$ for all $S\subseteq V$.

This follows by induction on the cardinality of $S$.
Each vertex is incident with at least three faces, so $f(S)\ge 3$.
Thus we may assume $|S|\ge 4$. There are two cases.
 If every face in $F$ is incident with $S$, then since a triangulation
on $s$ vertices has at least $2s-4$ faces, it follows that $f(S) \ge 2|S|-4$. This is
at least $|S|$ since $|S|\ge 4$. On the other hand, suppose there is a face in $F$
not containing a vertex of $S$. {We now shade all faces in $F$ that contain vertices of $S$. There must be an edge $xy$ that separates a shaded and unshaded face. Let $v$ be the third vertex of  the shaded face different from $x$ and $y$. We note that such a shaded face contains exactly one vertex of $S$, namely the vertex $v$.} By the induction hypothesis, it holds that $f(S-\{v\}) \ge |S|-1$, and $S$ is
incident with at least one more face, so that $f(S)\ge |S|$, as required.
\end{Proof}

Note that $T(G)$ has minimum degree $3$, and thus $\td( T(G) )\le 3$.

\begin{theorem}
If $G$ is a (not necessarily planar) triangulation of order at least $4$, then $\td( T(G) ) = 3$
if and only if $G$ is $3$-colorable.
\end{theorem}
\begin{Proof}
Suppose that $\td( T(G) ) = 3$. Let $D_1, D_2, D_3$ be a partition of $V( T(G) )$
into three vertex-disjoint dominating sets. Color the vertices of $G$ with $\{1,2,3\}$
by coloring each vertex $v$ with the index of the $D_i$ such that $v\in D_i$.
Each vertex of $T(G)$ not in $G$ has degree exactly~$3$, and those three
neighbors are in a different $D_i$. That is, every face of $G$ is
properly colored by our coloring. Since every edge of $G$ is in a face, it follows
that we have a proper $3$-coloring. That is, $G$ is $3$-colorable.

Conversely, suppose that $G$ has a proper $3$-coloring. Then, each color class $C$ must
appear in every face. Thus, every vertex of $T(G)$ not in $G$  has a neighbor in $C$,
as does every vertex of $G$ that has a color different from $C$.
It therefore suffices to show that for each vertex $v$ of $G$, we can find a
new vertex to color with the same color as~$v$. That is, we need a matching from
the vertices of $G$ to the faces of $G$ that saturates the vertices of~$G$.
Such a matching exists by Lemma~\ref{lem:match}.
\end{Proof}

We consider planar triangulations next.

\begin{theorem} \label{t:tGplanarLower}
If $G$ is a planar triangulation, then $\td( T(G) ) \ge 2$.
\end{theorem}
\begin{Proof}
By the Four Color Theorem, or for example~\cite{Penaud75planaires},
one can $2$-color the vertices of the base triangulation $G$ such that no
face is monochromatic. Each color class $C$ in such a $2$-coloring of $G$
totally dominates all vertices of $T(G)$ not in $G$, as well as all vertices of
$G$ whose color is different from $C$. Further, by Lemma~\ref{lem:match}, we
can match each vertex $v$ of $C$ with a new vertex $f_v$ of $T(G)$ not in $G$
and color $f_v$ with color $C$. In this way, every
vertex of $C$ has a neighbor in $C$. Thus, both color classes form a total
dominating set of~$T(G)$.
\end{Proof}

By Theorem~\ref{t:tGplanarLower}, it follows that $\FTD( T(G) ) \ge 2$
for all planar triangulations $G$. We know of four examples of equality.

\begin{lemma} \label{lem:ftd=2}
Let $G$ be one of the following four planar triangulations: $K_4$, the icosahedron,
or one of the two graphs drawn in Figure~\ref{f:kleetope2}. Then $\FTD( T(G) ) = 2$.
\end{lemma}

\begin{myFigure}
\includegraphics{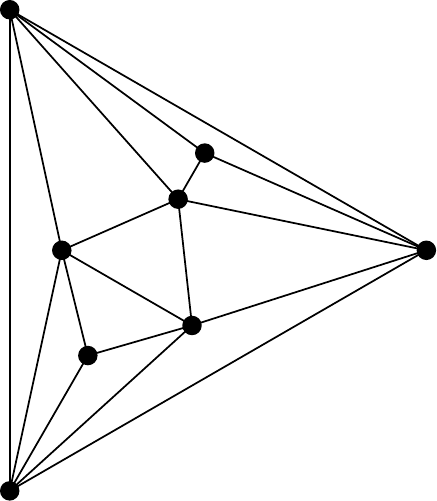} \qquad\qquad
\includegraphics{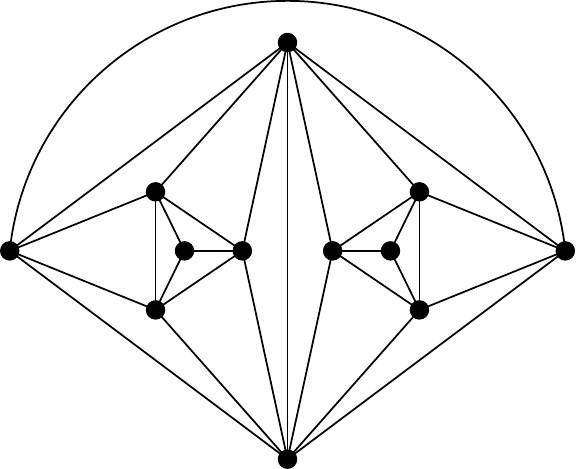}
\\
\vskip 0.2cm
\caption{Graphs whose Kleetopes have fractional total domatic number~$2$}
\label{f:kleetope2}
\end{myFigure}

\begin{Proof}
By Theorem~\ref{t:tGplanarLower}, we have $\td ( T(G) ) \ge 2$. So it suffices to show that $\FTD( T(G) ) \le 2$.

{Suppose firstly that $G = K_4$ and consider the Kleetope $T(G)$. In order to totally dominate the new vertices added in each face of the base graph $G$, all total dominating sets of $T(G)$ contain at least two vertices of $G$. Hence, if $\cF$ is a thoroughly distributed family of $T(G)$, then, by averaging over the four vertices in the base graph $G$, we have that $r_{_{\cF}}  \ge 2|\cF| / 4$, implying that the effective ratio $|\cF|/r_{_{\cF}}$ is at most~$2$ and hence $\FTD(T(G)) \le 2$.}

Similarly, with a bit more effort (or a computer) it can be shown that
for the graph of order~$8$ in Figure~\ref{f:kleetope2}, and for the icosahedron,
again every total dominating set of $T(G)$ contains at least half
the vertices of the base graph $G$, so that
it follows
that every thoroughly distributed family has effective ratio at most $2$.

Such a property is not true for the graph of order $12$ in Figure~\ref{f:kleetope2}. 
Instead we used a computer to calculate the value of the parameter
$\FTD(T(G))$.
\end{Proof}

We briefly consider the case that $G$ is a quadrangulation.

\begin{lemma}
If $G$ is a planar quadrangulation, then $\td ( T(G) ) \le 3$.
\end{lemma}
\begin{Proof}
Suppose that $T(G)$ has four disjoint total dominating sets.
As before, such a total domatic partition would,
when restricted to $G$, be a proper $4$-coloring of $G$. Each vertex
of $G$ thus needs a neighbor of the same color, which must be one of
the faces of $G$. But a facial vertex of $T(G)$ has four neighbors of different
colors, so it can be the
same color as only one neighbor. Thus we need a matching from $V(G)$ to $F(G)$ that saturates
$V(G)$, but in a planar quadrangulation, $|V(G)| > |F(G)|$
and therefore there is no such matching.
\end{Proof}

Our second construction is as follows.
For a triangulation $G$, the triangulation
$U(G)$ is obtained from the graph $G$
by adding in each face $f$ of $G$ a new triangle $\{f_1,f_2,f_3\}$ each vertex of which has
a different pair of neighbors on the boundary of the face. See Figure~\ref{f:triangulation}, where the white
vertices are new. As before, we will call $G$ the base graph.

 \begin{myFigure}
\includegraphics{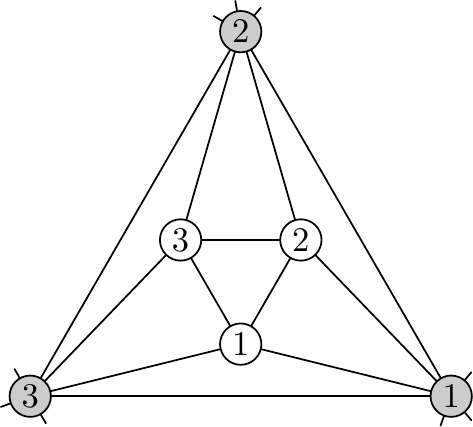} 
\\
\vskip 0.2cm
\caption{Construction of a  triangulation $U(G)$}
\label{f:triangulation}
\end{myFigure}

\begin{lemma} \label{lem:sg=3}
If $G$ is a planar triangulation, then $\td ( U(G) ) = \FTD( U(G) ) = 3$.
\end{lemma}
\begin{Proof}
By the result of Matheson and Tarjan~\cite{MaTa96}, the base graph $G$
has three disjoint dominating sets, say $D_1$, $D_2$, and $D_3$. These
can be extended into three disjoint total dominating sets of $U(G)$ as follows.
For each face $f$ of $G$, add one vertex $f_i$ of $U(G)$ not in $G$ to each $D_i$ such that
$f_i$ has a neighbor in $D_i$. Thus, $\td( U(G)) \ge 3$.

On the other hand, consider a face $f$ of $G$. It is easy to verify that
every total dominating set of $U(G)$ must contain at least two of
the vertices from the $6$-tuple consisting of the boundary of $f$ and the new triangle inside $f$.
Thus, by
averaging, the effective ratio of any thoroughly distributed family of $U(G)$ is at most~$3$.
\end{Proof}

\subsection{Lower Bounds for Triangulations}

We now turn to lower bounds for triangulations. Unfortunately, we mostly
have only open questions.
Since every planar triangulation is a triangulated disc, Theorem~\ref{t:disc}
implies that every planar triangulation $G$ satisfies $\FTD(G) \ge \frac{3}{2}$.
We believe this lower bound can be improved significantly and pose the
following conjecture.

\begin{conjecture}
\label{conj:triangulationHasTwo}
If $G$ is a planar triangulation of order at least~$4$, then $\td(G) \ge 2$.
\end{conjecture}

We can establish the conjecture for a few cases.

\begin{observation}
If $G$ is a planar triangulation where
every vertex has odd degree, then $\td(G) \ge 2$.
\end{observation}
\begin{Proof}
This follows from the Four Color Theorem.
If a vertex $v$ has
odd degree in a triangulation, its neighborhood contains an odd cycle.
Thus in a proper coloring, $v$ has neighbors of three different colors.
Thus, the union of two color classes is a total dominating set.
\end{Proof}

Now, it is not true in general that if one has a proper $4$-coloring of a
planar triangulation, then one can just combine two color-classes to obtain a
total dominating set. However, computer search suggests that:

\begin{conjecture} \label{conj:fourColorExtension}
Every planar triangulation with at least four vertices has
a proper $4$-coloring $(C_1,C_2,C_3,C_4)$ such that $C_1 \cup C_2$ and $C_3
\cup C_4$ are total dominating sets.
\end{conjecture}

Equivalently, {$V(G)$ can be partitioned into two total dominating sets both of which induce a bipartite subgraph of $G$.} 

Another case where Conjecture~\ref{conj:triangulationHasTwo}
(and Conjecture~\ref{conj:fourColorExtension}) holds is where the dual is hamiltonian.

\begin{observation}
If $G$ is a planar triangulation and the dual of $G$ is hamiltonian,
then $\td(G) \ge 2$.
\end{observation}
\begin{Proof}
This uses a standard idea (as for example used by Tait).
Consider a hamilton cycle of the dual of $G$ as a closed curve $\cC$ in the plane.
Color vertices of $G$ red if they are inside the curve $\cC$, and blue if they are
outside the curve. We claim that every vertex $v$ has both a red and a blue
neighbor.
For, consider the cycle through $N(v)$
as a closed curve~$\cD$. Then, the curve $\cC$ must cut $\cD$ (in at least
two places). Say it cuts the edge $xy$. Then $x$ and $y$ have different colors.
\end{Proof}

Now, the dual of a planar triangulation is a $3$-connected cubic planar graph.
The smallest such graphs that are not hamiltonian have 38 vertices. (See~\cite{HoMc88planar}.)
It can be checked that the duals of these nonhamiltonian graphs do have two
disjoint total dominating sets. By the technique of Penaud~\cite{Penaud75planaires}, to prove
Conjecture~\ref{conj:triangulationHasTwo}
it would be enough to show that: every $3$-connected cubic planar graph has a
$2$-factor that does not include a facial cycle.

Now, there are planar triangulations that have total domatic number $2$.
For example, if we take $T(G)$ for any planar triangulation with chromatic
number $4$. But there are also planar triangulations that have
fractional total domatic number~$2$, though the four listed in Lemma~\ref{lem:ftd=2}
are the only we know of. This shows that Conjecture~\ref{conj:triangulationHasTwo} if true is sharp.

We remark that the Kleetope construction shows that
Conjecture~\ref{conj:triangulationHasTwo} does not carry over to graphs embeddable on
the torus or Klein bottle, as the Kleetope, $T(K_5)$, of $K_5$ does not have
two disjoint total dominating sets. That is:

\begin{observation}
There exists toroidal triangulations $G$ satisfying $\td(G) = 1$.
\end{observation}

If one imposes larger minimum degree, it appears even more can be said.
We pose the following conjecture.

\begin{conjecture}
If $G$ is a planar triangulation with minimum degree at least~$4$, then
$\td( G ) \ge 3$.
\end{conjecture}

\noindent
If true, this conjecture is sharp by the graphs $U(G)$ defined earlier. See Lemma~\ref{lem:sg=3}.

We conclude this section with even more speculation.
Perhaps it is true that every triangular disc with minimum degree
at least~$3$ has two disjoint total dominating sets.
It is not true that  every triangular disc with minimum degree
at least~$4$ has three disjoint total dominating sets: the icosahedron
minus a vertex is an example, and there is an example of order 10.
But maybe there are only finitely many exceptions.

\section{Cubic Graphs}

In this section we consider the fractional total domatic number of cubic graphs.
The question of which cubic graphs have two disjoint total dominating sets
is well studied. For example, it is known that
the Heawood graph is the smallest example without two disjoint total dominating
sets. (For more information see for example McCuaig~\cite{McCuaig00even},
or Gropp~\cite{Gropp97blocking}.)
In contrast, Thomassen~\cite{Thomassen92} showed that, for $r\ge 4$, every $r$-regular graph
has two disjoint total dominating sets.

It is known \cite{HY13hypergraph} that a connected $3$-regular $3$-uniform hypergraph is almost
$2$-colorable, in that there is a $2$-coloring that $2$-colors all but one
specified hyperedge. Similarly, a connected cubic graph has a $2$-coloring
such that each vertex, except possibly two, has both colors in its neighborhood.
This enables us to provide a lower bound on the fractional
parameters.

\begin{theorem}
If $H$ is a connected $3$-regular $3$-uniform hypergraph on $n$ vertices, then $\FdT(H) \ge 2n/(n+1)$.
\end{theorem}
\begin{Proof}
Let $e$ be
any hyperedge containing~$v$ and let $H_v$ be the subhypergraph with~$e$ removed from $H$.
By the result of~\cite{HY13hypergraph}, the hypergraph $H_v$ has a $2$-coloring,
say $(A_v, B_v)$. If we add $v$ to
whichever of $A_v$ and $B_v$ doesn't already contain it, then both $A_v$ and~$B_v$
are transversals of $H$.

Consider the collection $\cF = \{\, A_v : v\in V(H) \,\} \cup \{\, B_v : v \in V(H) \,\}$.
Each vertex $v$ is in at most $n+1$ of these sets. So
$\cF$ is a transversal family with effective transversal-ratio $2n/(n+1)$,
as required.
\end{Proof}

By the connection with the ONH of $G$, it follows that
a connected cubic graph $G$ of order $n$ has $\FTD(G) \ge 2-o(1)$.
That is, cubic graphs ``almost'' have two disjoint total dominating sets.
However, we are unable to find examples of equality in the above theorem
or cubic graphs with fractional total domatic number less than $2$.
So we conjecture:

\begin{conjecture} \label{conj:cubic}
If $G$ is a connected cubic graph, then $\FTD(G) \ge 2$.
Indeed, $G$ has a thoroughly distributed family of four sets such that
every vertex is in at most~two of these.
\end{conjecture}

There is computer support for this.
Indeed, the computer suggests that $\FTD(G) > 2$
for all connected cubic graphs $G$ except those that have $\gt(G) = n/2$ (which have
$\td(G) = \FTD(G) = 2$). (See~\cite{HeYe08} for the description of such graphs.)
Or more generally, we ask:

\begin{question}
Does every $3$-uniform $3$-regular hypergraph have
four transversals with each vertex in at most two?
\end{question}

\section{Open Questions}

A lot remains to be ascertained. The two questions most frustrating to
us are Conjecture~\ref{conj:triangulationHasTwo}, that every planar triangulation has two
disjoint total dominating sets, and Conjecture~\ref{conj:cubic}, that every
cubic graph has fractional total domatic number at least $2$. We
also wonder about the relationship between the total domatic
number and its fractional counterpart. For example, how large
can $\FTD(G)$ be if $\td(G)=1$? The largest we know of
is $7/3$ from the Heawood graph.


\begin{thebibliography}{99}

\bibitem{AELTWarxiv}
W.~Abbas, M.~Egerstedt, C.-H. Liu, R.~Thomas, and P.~Whalen,
\emph{Deploying robots with two sensors in ${K}_{1,6}$-free graphs},
J. Graph Theory  \textbf{82}  (2016),  236--252.


\bibitem{ASV12regular}
H. Aram,  S.M. Sheikholeslami and L. Volkmann,
\emph{On the total domatic number of regular graphs},
Trans. Comb.  \textbf{1}  (2012), 45--51.

\bibitem{Berge79survey}
C.~Berge,
\emph{Packing problems and hypergraph theory: a survey},
Ann. Discrete Math. \textbf{4} (1979), 3--37.

\bibitem{BuHeTu12}
C.~Bujt\'{a}s, M.A. Henning, and Zs. Tuza,
\emph{Transversals and domination in uniform hypergraphs},
European J. Combin. \textbf{33} (2012), 62--71.

\bibitem{BuHeTuYe14}
C.~Bujt\'{a}s, M.A. Henning, Zs. Tuza, and A.~Yeo,
\emph{Total transversals and total domination in uniform hypergraphs},
Electronic J. Combin. \textbf{21(2)} (2014), \#P2.24.

\bibitem{CKTV15coupon}
B.~Chen, J.H. Kim, M.~Tait, and J.~Verstraete,
\emph{On coupon colorings of graphs},
Discrete Applied Mathematics \textbf{193} (2015), 94--101.

\bibitem{CDH80total}
E.J. Cockayne, R.M. Dawes, and S.T. Hedetniemi,
\emph{Total domination in graphs},
Networks \textbf{10} (1980), 211--219.


\bibitem{DHJ16outerplanar}
M.~Dorfling, J.H. Hattingh, and E.~Jonck,
\emph{Total domination in maximal outerplanar graphs {II}},
Discrete Math. \textbf{339} (2016), 1180--1188.

\bibitem{FH08clawfreecubic}
O.~Favaron and M.A. Henning,
\emph{Bounds on total domination in claw-free cubic graphs},
Discrete Math. \textbf{308} (2008), 3491--3507.

\bibitem{FYK00configurations}
S.~Fujita, M.~Yamashita, and T.~Kameda,
\emph{A study on {$r$}-configurations---a resource assignment problem on graphs},
SIAM J. Discrete Math. \textbf{13} (2000), 227--254 (electronic).

\bibitem{Furedi81hypergraphs}
Z.~F{\"u}redi,
\emph{Maximum degree and fractional matchings in uniform hypergraphs},
Combinatorica \textbf{1} (1981), 155--162.

\bibitem{Gropp97blocking}
H.~Gropp,
\emph{Blocking set free configurations and their relations to digraphs and hypergraphs},
Discrete Math. \textbf{165/166} (1997), 359--370.

\bibitem{HeTe98}
P.~Heggernes and J.~A. Telle,
\emph{Partitioning graphs into generalized dominating sets},
Nordic J. Comput. \textbf{5} (1998), 128--142.

\bibitem{HS08disjoint}
M.A. Henning and J.~Southey,
\emph{A note on graphs with disjoint dominating and total dominating sets},
Ars Combin. \textbf{89} (2008), 159--162.

\bibitem{HeYe08}
M.A. Henning and A.~Yeo,
\emph{Hypergraphs with large transversal number and with edge sizes at least three},
J. Graph Theory \textbf{59} (2008)  326--348.

\bibitem{HY13hypergraph}
M.A. Henning and A.~Yeo,
\emph{2-colorings in {$k$}-regular {$k$}-uniform hypergraphs},
European J. Combin. \textbf{34} (2013), 1192--1202.

\bibitem{HoMc88planar}
D.A. Holton and B.D. McKay,
The smallest non-Hamiltonian $3$-connected cubic planar graphs have $38$ vertices.
J. Combin. Theory Ser. B  \textbf{45}  (1988),  305--319.

\bibitem{KoWo01panchromatic}
A.V. Kostochka and D.R. Woodall,
\emph{Density conditions for panchromatic colourings of hypergraphs},
Combinatorica  \textbf{21}  (2001), 515--541.

\bibitem{Lovasz70konig}
L. Lov{\'a}sz, \emph{A generalization of {K}\"onig's theorem}, Acta Math. Acad. Sci. Hungar. \textbf{21} (1970), 443--446.

\bibitem{MaTa96}
L.R. Matheson and R.E. Tarjan,
\emph{Dominating sets in planar graphs},
European J. Combin. \textbf{17} (1996), 565--568.

\bibitem{McCuaig00even}
W.~McCuaig,
\emph{Even dicycles},
J. Graph Theory \textbf{35} (2000), 46--68.

\bibitem{Penaud75planaires}
J.G. Penaud,
\emph{Une propri\'et\'e de bicoloration des hypergraphes planaires},
Cahiers Centre \'Etudes Recherche Op\'er. \textbf{17} (1975), 345--349.

\bibitem{ScUl97fractional}
E.R. Scheinerman, and D.H. Ullman,
\emph{Fractional graph theory:  A rational approach to the theory of graphs},
Wiley, 1997.

\bibitem{Suomela06sleep}
J.~Suomela,
\emph{Locality helps sleep scheduling},
Working Notes of the
  Workshop on World-Sensor-Web: Mobile Device-Centric Sensory Networks and
  Applications (WSW, Boulder, CO, USA, October 2006), p.~41--44.

\bibitem{ThYe07}
S.~Thomass\'{e} and A.~Yeo,
\emph{Total domination of graphs and small transversals of hypergraphs},
Combinatorica \textbf{27} (2007), 473--48.

\bibitem{Thomassen92}
C. Thomassen,
\emph{The even cycle problem for directed graphs},
 J. Amer. Math. Soc.  \textbf{5} (1992),  217--229.

\bibitem{Zelinka89total}
B.~Zelinka, \emph{Total domatic number and degrees of vertices of a graph},
  Math. Slovaca \textbf{39} (1989), 7--11.

\end{thebibliography}
\end{document}